\documentclass{amsart}
\usepackage{graphicx}
\title{Properly embedded and immersed minimal surfaces in the Heisenberg group}
\author{ Jih-Hsin Cheng }
\address[Cheng]
{Institute of Mathematics,
Academia Sinica, \newline%
\indent Nankang,
Taipei, Taiwan, 11529,
R.O.C.}
\email[]{cheng@math.sinica.edu.tw}%
\author{Jenn-Fang Hwang}
\address[Hwang]
{Institute of Mathematics,
Academia Sinica, \newline%
\indent Nankang,
Taipei, Taiwan, 11529,
R.O.C.}
\email[]{majfh@ccvax.sinica.edu.tw}%
\date{}

\begin{document}
\maketitle
\centerline{\bf Abstract}
We study properly embedded and immersed p(pseudohermitian)-minimal surfaces in the 3-dimensional Heisenberg group. From the recent work of Cheng, Hwang, Malchiodi, and Yang, we learn that such surfaces must be ruled surfaces. There are two types of such surfaces: band type and annulus type according to their topology. We give an explicit expression for these surfaces. Among band types there is a class of properly embedded p-minimal surfaces of so called helicoid type. We classify all the helicoid type p-minimal surfaces. This class of p-minimal surfaces includes all the entire p-minimal graphs (except contact planes) over any plane. Moreover, we give a necessary and sufficient condition for such a p-minimal surface to have no singular points. For general complete immersed p-minimal surfaces, we prove a half space theorem and give a criterion for the properness.
 
\bigskip

\section{{\bf Introduction and statement of the results}}

In [CHMY] we developed a surface theory in pseudohermitian geometry. In particular we defined the notion of p(pseudohermitian)-minimal surface. The equation for a graph $(x,y,u(x,y))$ in $R^3$ to be p-minimal reads

\begin{align}
&div(\frac{\nabla u+\vec{F}}{|\nabla u+\vec{F}|})=0&\nonumber
\end{align}

\noindent where $\vec{F}$ denotes the plane vector (field) $(-y,x)$ (associated with the standard contact structure of $R^3$). As a differential equation, the above p-minimal surface equation is degenerate (hyperbolic and elliptic).
By analyzing the singular set (where $\nabla u+\vec{F}=0$), we solved the analogue of the Bernstein problem in the 3-dimensional Heisenberg group $H_1$ (identified with $R^3$ as a set). Namely we showed that two known families of examples ([Pau]) are the only entire $C^2$ smooth p-minimal graphs over the $xy$-plane. In this paper we want to study general $C^2$ smooth p-minimal surfaces properly embedded in $H_1$ (we will often omit "$C^2$ smooth" hereafter). 

First let us recall the history of (Riemannian or Euclidean) minimal surface theory in $R^3$ briefly (e.g., [Oss], [Nit], [DHKW], [Fan], etc.). Until 1980 only a few complete embedded minimal surfaces had been found like, for instance, the helicoid, the catenoid, Scherk's surface, Riemann's examples, and so on. Among these surfaces the catenoid and the helicoid were the only known complete embedded minimal surfaces with finite genus in $R^3$ at the turn of 1980. It was only in 1982 that Costa ([Cos]) discovered a new complete minimal surface immersed in $R^3$ of genus one. This example was shown to be embedded soon after it was discovered ([HM1]). Subsequently various complete embedded minimal surfaces were discovered ([HM1], [HKW], etc.; see also [DHKW]). In 1997 Collin ([Col]) proved that the catenoid is the only annulus type, properly embedded minimal surface. Around the year 2000, it was announced that a properly embedded, simply connected minimal surface in $R^3$ is either a plane or a helicoid ([MR]). 

In contrast to the Riemannian case, a properly embedded p-minimal surface in $R^3$, identified with $H_1$, must be a ruled surface with the rulings generated by Legendrian lines ([CHMY]). A general ruled surface satisfies a 3rd order equation (see page 225 in [Mon] or the last paragraph of Section 4 in [CHMY] for a brief explanation). Requiring the rulings lying in contact planes (this is what "Legendrian" means) restricts such a surface to satisfy a 2nd order equation. It is not hard to see that a connected, properly embedded ruled surface must be homeomorphic to either $R^2$ (band type) or $R^{1}{\times}S^{1}$ (annulus type). So there are no properly embedded p-minimal surfaces of positive genus type in $R^3$. 

A band type, properly embedded p-minimal surface is called of helicoid type if its Legendrian rulings are lying in the parallel planes and the union of all such parallel planes is the whole space $H_1$ (if we remove the latter restriction, we call such a surface of helicoid type in the weak sense, cf. Example 4.1). This includes the class of entire p-minimal graphs (excluding the contact planes) over any plane. To classify all such entire $C^2$ smooth p-minimal graphs (by "entire" we mean "defined on the whole plane"), we first extend the class of such graphs to the class of helicoid type. Then we find out all the helicoid type p-minimal surfaces (see Theorem B). Let us explain this idea in detail.
 
Let $\Sigma$ be an entire $C^2$ smooth p-minimal graph over a plane $P$ in $H_1$, identified with the Euclidean space $R^3$. As mentioned above $\Sigma$ is a classical ruled surface with the rulings generated by Legendrian lines, called characteristic lines. If two such lines intersect, we can show that $\Sigma$ must be a contact plane past the intersection point (see Proposition 2.1). Otherwise we can project all characteristic lines to a family of parallel lines $\Gamma_t$ on $P$ (note that $\Sigma$ is a graph over $P$). For each $\Gamma_t$ we can find a unique plane $P_t$ perpendicular to $P$ and containing $\Gamma_t$. It follows that all $P_t$'s are parallel to each other. Thus we can extend the problem of finding entire $C^2$ smooth p-minimal graphs (excluding the contact planes) to the following one:

\medskip

{\bf Problem A.} {\it Suppose $\{ P_t$: $t\in R \}$ is a family of disjoint parallel planes in $H_1$ such that $\bigcup_{t\in R}P_{t}=H_{1}$. Find all $C^2$ smooth p-minimal surfaces $\Sigma$ such that ${\tilde{\Gamma}}_{t}\equiv\Sigma\cap P_{t}$ is a characteristic (whole straight) line of $\Sigma$ for each $t$. Namely (see Section 2) find all $C^2$ smooth, helicoid type p-minimal surfaces in $H_1$.}

\medskip

Take a plane $P$ perpendicular to $P_t$ for all $t$. Some characteristic lines ${\tilde{\Gamma}}_{t}$ may be perpendicular to $P$. So a solution surface to Problem A may not be a graph over a certain plane. 

In Section 2 we start with a study of general immersed and properly embedded p-minimal surfaces in $H_1$.  Define $X$ by (2.1), (2.2), and (2.3), describing a $C^2$ smooth p-minimal surface in $H_1$. We obtain a necessary and sufficient condition for $X$ being embedded (implied by $X$ being immersed, injective, and proper). We also obtain a criterion for a point to be singular (see the definition in [CHMY] and the generalized definition in the proof of Theorem A (b)).

\bigskip

\noindent{\bf Theorem A.} {\it (a) $X$ is immersed if and only if for any $(s,t)$ either (2.8a) or (2.8b) fails to hold. (b) $X$ is singular at $(s,t)$ if and only if (2.9) holds. (c) If we take $\xi =0$ in $(2.7)$, then $X$ is injective if and only if for ${t_1}\neq{t_2}$ $\hat{z}_{2}-\hat{z}_{1}$ in (2.15) is not 0 in the case that $\Gamma_{t_1}\cap\Gamma_{t_2}$ consists of exactly one point while $\gamma_{1}\neq\gamma_{2}$ in the case that $\Gamma_{t_1}=\Gamma_{t_2}$.}
  
\bigskip

Note that we do not lose generality by assuming $\xi =0$ in Theorem A (c). We give a geometric interpretation (see (2.16)) for the formula (2.15) and a direct application (see (2.18)). 

In Section 3 we deal with the helicoid type p-minimal surfaces. We observe that
the following are two families of solutions to Problem A:

\begin{align}
&u=-abx^{2}+(a^{2}-b^{2})xy+aby^{2}+g(-bx+ay)& \tag{1.1} \\
&(a,b \mbox{ being real constants such that }a^{2}+b^{2}=1 \mbox{ and }  g\in C^{2});&\nonumber\\
&(x-x_{0})cos{\theta}(t)+(y-y_{0})sin{\theta}(t)=0& \tag{1.2} \\
&\mbox{where }t=z-y_{0}x+x_{0}y,\ \theta\in C^{2}(R),\mbox{ and
}x_{0},y_{0}\mbox{ are real constants}&\nonumber
\end{align}

\noindent (see also [Pau] for (1.1)). The following result resolved Problem A.

\bigskip

\noindent{\bf Theorem B.} {\it Any $C^2$ smooth, helicoid type p-minimal surface in $H_1$ can be expressed in the form of either (1.1) or (1.2).}

\bigskip

We also give a criterion for (1.2) to have no singular points.

\bigskip

\noindent{\bf Theorem C.} {\it (1.2) has no singular points if and only if $\theta^{\prime}(t)\geq{0}$ for all $t\in{R}$.}

\bigskip

Note that (1.1) has no isolated singular points, but does have one (connected) singular curve ([CHMY]). By taking ${\theta}(t)=-cot^{-1}(t)$ and $x_{0}=y_{0}=0$ in (1.2), we obtain $y=xz$. This is an entire p-minimal graph over the $xz$-plane having no singular points. We remark that the vertical planes $ax+by+c=0$ have been the only known entire p-minimal graphs having no singular points before ([GP1]). After this paper was completed, we were informed that another example of entire p-minimal graph having no singular points and the result about entire p-minimal graphs over any plane (the graph case of our Theorem B) are also obtained in a new preprint ([GP2]). 

In Section 4, we discuss the properness (Proposition 4.1) and non-helicoid type p-minimal surfaces. The first possibility is that $\bigcup_{t\in R}P_{t}\neq H_{1}$ where $P_{t}'s$ are the parallel planes containing characteristic lines. This can occur as shown in Example 4.1. However, if we confine the surface to the upper half space $R^{3}_{+}=\{ (x,y,z)\in R^3\, :\, z>0\}$, then this is not possible unless the surface is a contact plane (see Theorem D below and observe that the $xy$-plane is the contact plane past the origin). A contact plane $P$ in $H_{1}$ divides $H_{1}$ into two halfspaces $P_{\pm}$ ($H_{1}\backslash P=P_{+}\cup P_{-}$). Since $P$ is not perpendicular to the $xy$-plane, we can talk about the {\it upper halfspace} $P_{+}$ (containing "$(0,0,+\infty )$") and the {\it lower halfspace}  
$P_{-}$ (containing "$(0,0,-\infty )$"). We have the following halfspace theorem.  

\bigskip 

\noindent {\bf Theorem D.} {\it Suppose that $\Sigma$ is a $C^2$ smooth, connected, complete, immersed p-minimal surface in the halfspace $P_{+}$ (or the halfspace $P_{-}$). Then $\Sigma$ is a contact plane parallel to $P$.}

\bigskip

\noindent We remark that the halfspace theorem in the Riemannian minimal surface theory was proved by Hoffman and Meeks III ([HM2]). 


\bigskip

{\bf Acknowledgments.} We would like to thank Andrea Malchiodi and Paul Yang for many discussions in this newly developed research direction. In particular, Andrea pointed out to us the examples with a multivalued function $g$ in $(1.1)$ (Example 4.1 is such an example).   

\bigskip

\section{{\bf General p-minimal surfaces and proof of Theorem A}}

We learned from [CHMY] that a $C^2$ smooth p-minimal surface in $H_1$ must be a ruled surface in $R^3$ (identified with $H_1$). Consider a connected, $C^2$ smooth, complete p-minimal surface $\Sigma$ immersed in $H_1$ (by connectedness we mean the domain of the immersion is connected). First we observe the following fact.

\bigskip

\noindent {\bf Proposition 2.1.} {\it Suppose there are two characteristic lines on $\Sigma$ meeting at a point $p$ locally. Then $p$ is a singular point, the tangent space $T_{p}\Sigma$ is the contact plane past $p$, and ${\Sigma}=T_{p}\Sigma$.}

\bigskip    
 
\noindent By "meeting at a point locally", we exclude the possible intersection due to the "return" of characteristic lines.

\medskip

{\bf Proof of Proposition 2.1:}

Locally we can express $\Sigma$ as a graph $\{ (x,y,u(x,y))\}$ over the $xy$-plane. Let $\bar p$ be the projection of $p$ on the $xy$-plane. First by the definition of characteristic curves (see Section 1 in [CHMY]), $p$ ($\bar p$) must be a singular point. By Theorem 3.3 in [CHMY], either $\bar p$ is an isolated point in the singular set $S(u)$ or $S(u)$ is a $C^1$ smooth curve in a neighborhood of $\bar p$. In the latter case, only one characteristic line can meet $\bar p$ according to Corollary 3.6 in [CHMY]. So this case is impossible. In the former case, the union of all characteristic lines passing through $p$, the contact plane past $p$, forms the tangent space $T_{p}\Sigma$ by Lemma 4.3 in [CHMY]. By completeness $\Sigma$ must contain $T_{p}\Sigma$. By connectedness $\Sigma$ must coincide with $T_{p}\Sigma$. 

\begin{flushright}
Q.E.D.
\end{flushright} 

Now let us start with (4.9) in [CHMY]: describe a $C^2$ smooth p-minimal
surface $\Sigma$ in $H_1$ in parameters $s$ and $t$ as follows:

\begin{align}
&x=s(sin{\theta}(t))+{\alpha}(t)&\tag{2.1}\\
&y=-s(cos{\theta}(t))+{\beta}(t)& \tag{2.2}\\
&z=s[{\beta}(t)sin{\theta}(t)+{\alpha}(t)cos{\theta}(t)]+{\gamma}(t)& \tag{2.3}
\end{align}  

\noindent where ${\alpha}$, ${\beta}$, ${\gamma}$, and $\theta$ are all in $C^2$. Let $X(s,t)=
(x(s,t),y(s,t),z(s,t))$. 

\medskip

{\bf Proof of part (a) of Theorem A:}

$X$ being an immersion implies a certain relation among
${\alpha}$, ${\beta}$, ${\gamma}$, and $\theta$. It is easy to see that ${\partial_s}X = (sin{\theta},
-cos{\theta}, {\beta}sin{\theta}+{\alpha}cos{\theta})$ and ${\partial_t}X = s{\theta}^{\prime}
(cos{\theta}, sin{\theta}, {\beta}cos{\theta}-{\alpha}sin{\theta}) + s(0,0,{\beta}^{\prime}
sin{\theta}+{\alpha}^{\prime}cos{\theta}) + ({\alpha}^{\prime}, {\beta}^{\prime},
{\gamma}^{\prime})$. Then a straightforward computation gives the $x,y,z$ components
of the cross product ${\partial_s}X\times{\partial_t}X$:

\begin{align}
&-s{\beta}{\theta}^{\prime}-scos{\theta}({\beta}^{\prime}
sin{\theta}+{\alpha}^{\prime}cos{\theta})-{\beta}^{\prime}({\beta}sin{\theta}+{\alpha}cos{\theta})-{\gamma}^{\prime}cos{\theta},&\tag{2.4}\\
&s{\alpha}{\theta}^{\prime}-ssin{\theta}({\beta}^{\prime}
sin{\theta}+{\alpha}^{\prime}cos{\theta})+{\alpha}^{\prime}({\beta}sin{\theta}+{\alpha}cos{\theta})-{\gamma}^{\prime}sin{\theta},&\tag{2.5}\\
&s{\theta}^{\prime}+{\beta}^{\prime}
sin{\theta}+{\alpha}^{\prime}cos{\theta},&\tag{2.6}
\end{align}

\noindent respectively. Write 





\begin{align}
&({\alpha}(t),{\beta}(t)) = {\delta}(t)(cos{\theta}(t),sin{\theta}(t))+{\xi}(t)(sin{\theta}(t),-cos{\theta}(t))&\tag{2.7}
\end{align}


\noindent where ${\delta}(t) = {\alpha}(t)cos{\theta}(t) + {\beta}(t)sin{\theta}(t)$ and ${\xi}(t) = {\alpha}(t)sin{\theta}(t) - {\beta}(t)cos{\theta}(t)$. Substituting (2.7) into (2.1), (2.2), and (2.3) gives

\begin{align}
&x=(s+\xi (t))sin\theta (t)+\delta (t)cos\theta (t)&\tag{$2.1^{\prime}$}\\
&y=-(s+\xi (t))cos\theta (t)+\delta (t)sin\theta (t)&\tag{$2.2^{\prime}$}\\
&z=s\delta (t)+\gamma (t).&\tag{$2.3^{\prime}$}
\end{align}

\noindent We compute ${\partial_s}X = (sin\theta , -cos\theta , \delta )$ and ${\partial_t}X = (Acos\theta + Bsin\theta , Asin\theta - Bcos\theta , s{\delta}^{\prime} + {\gamma}^{\prime})$ where $A=(s+\xi ){\theta}^{\prime} + {\delta}^{\prime}$ and $B={\xi}^{\prime} - {\delta}{\theta}^{\prime}$. Then the $x,y,z$ components
of the cross product ${\partial_s}X\times{\partial_t}X$ read

\begin{align}
&-cos\theta (s{\delta}^{\prime} + {\gamma}^{\prime}-\delta B) - A\delta sin\theta , &\tag{$2.4^{\prime}$}\\
&-sin\theta (s{\delta}^{\prime} + {\gamma}^{\prime}-\delta B) + A\delta cos\theta ,
&\tag{$2.5^{\prime}$}\\
&A,&\tag{$2.6^{\prime}$}
\end{align}

\noindent respectively. Now we can deduce the condition for a point where $X$ is not immersed (${\partial_s}X\times{\partial_t}X=0$) as follows:

\begin{align}
&(s+\xi (t)){\theta}^{\prime}(t) + {\delta}^{\prime}(t)=0&\tag{2.8a}\\
&s{\delta}^{\prime}(t) + {\gamma}^{\prime}(t)-\delta (t)({\xi}^{\prime}(t) - {\delta}(t){\theta}^{\prime}(t))=0.&\tag{2.8b}
\end{align}

\noindent So $X$ is an immersion if and only if for any $(s,t)$ either (2.8a) or (2.8b) fails to hold. We have proved part (a) of Theorem A.

\begin{flushright}
Q.E.D.
\end{flushright} 

\medskip

{\bf Proof of part (b) of Theorem A:}

First observe that a singular point is where ${\partial_s}X\times{\partial_t}X$ is parallel to the vector $(-y, x, 1)$ (the vector dual to the contact form $dz + xdy - ydx$) (we may take this property as the generalized definition of a singular point where $X$ may not even be an immersion). It follows that $s{\delta}^{\prime} + {\gamma}^{\prime}-\delta B = -(s + \xi )A$ at a singular point by comparing $(2.4^{\prime})$, $(2.5^{\prime})$, $(2.6^{\prime})$ with $(2.2^{\prime})$, $(2.1^{\prime})$. Substitution of $A$ and $B$ gives

\begin{align}
&[(s+\xi )^{2}+{\delta}^{2}]{\theta}^{\prime}+(2s+\xi ){\delta}^{\prime}+{\gamma}^{\prime}-\delta {\xi}^{\prime}=0.&\tag{2.9}
\end{align} 

\begin{flushright}
Q.E.D.
\end{flushright} 

\medskip

Next we will discuss the injectivity of $X$. Eliminating $s$ in (2.1), (2.2), and (2.3), we obtain

\begin{align}
&(x-{\alpha}(t))cos{\theta}(t)+(y-{\beta}(t))sin{\theta}(t)=0&\tag{2.10}\\
&z={\beta}(t)x-{\alpha}(t)y+{\gamma}(t).&\tag{2.11}  
\end{align} 

\noindent For a fixed $t$, (2.10) is the equation for the $xy$-plane projection $\Gamma_{t}$ of the characteristic line $\tilde{\Gamma}_t$ while (2.11) together with (2.10) describes $\tilde{\Gamma}_t$. Suppose $\Gamma_{t_1}\cap\Gamma_{t_2}$ consists of exactly one point $(\hat{x},\hat{y})$ for ${t_1}\neq{t_2}$. Write ${\alpha}(t_{i}), {\beta}(t_{i}), {\gamma}(t_{i}), {\theta}(t_{i})$ as ${\alpha}_{i}, {\beta}_{i}, {\gamma}_{i}, {\theta}_{i}$, respectively for $i=1,2$. Then from (2.10) we can easily deduce

\begin{align}
&\hat{x}=\frac{{\alpha}_{1}cos\theta_{1}sin\theta_{2}-\alpha_{2}cos\theta_{2}sin\theta_{1}-(\beta_{2}-\beta_{1})sin\theta_{1}sin\theta_{2}}{sin(\theta_{2}-\theta_{1})},&\tag{2.12}\\
&\hat{y}=\frac{\beta_{2}sin\theta_{2}cos\theta_{1}-\beta_{1}sin\theta_{1}cos\theta_{2}+(\alpha_{2}-{\alpha}_{1})cos\theta_{1}cos\theta_{2}}{sin(\theta_{2}-\theta_{1})}.&\tag{2.13}
\end{align}

\noindent Let $\hat{z}_{i}={\beta}_{i}\hat{x}-{\alpha}_{i}\hat{y}+{\gamma}_{i}$ denote the $z$-coordinate of the point in $\tilde{\Gamma}_{t_i}$, projected down to $(\hat{x},\hat{y})$ in the $xy$-plane by (2.11) for $i=1,2$. We can then compute the difference of $\hat{z}_{1}$ and $\hat{z}_{2}$ by the substitution of (2.12) and (2.13):

\begin{align}
\hat{z}_{2}-\hat{z}_{1}&=\frac{-[(\alpha_{2}-{\alpha}_{1})cos\theta_{2}+
(\beta_{2}-\beta_{1})sin\theta_{2}][(\alpha_{2}-{\alpha}_{1})cos\theta_{1}+(\beta_{2}-\beta_{1})sin\theta_{1}]}{sin(\theta_{2}-\theta_{1})}&\tag{2.14}\\
&-(\alpha_{2}\beta_{1}-{\alpha}_{1}\beta_{2})+ {\gamma}_{2}-{\gamma}_{1}.&\nonumber
\end{align}

\noindent Let ${\psi}=\theta_{2}-\theta_{1}$. Denote the points $(\hat{x},\hat{y})$, $({\alpha}_{1},\beta_{1})$, $({\alpha}_{2},\beta_{2})$, $(0,0)$ by $P, Q_{1}, Q_{2}, O$, respectively. Assume that $PQ_{1}OQ_{2}P$ is in the order of the counterclockwise direction as shown in Figure 1. 

\includegraphics[width=8.0cm]{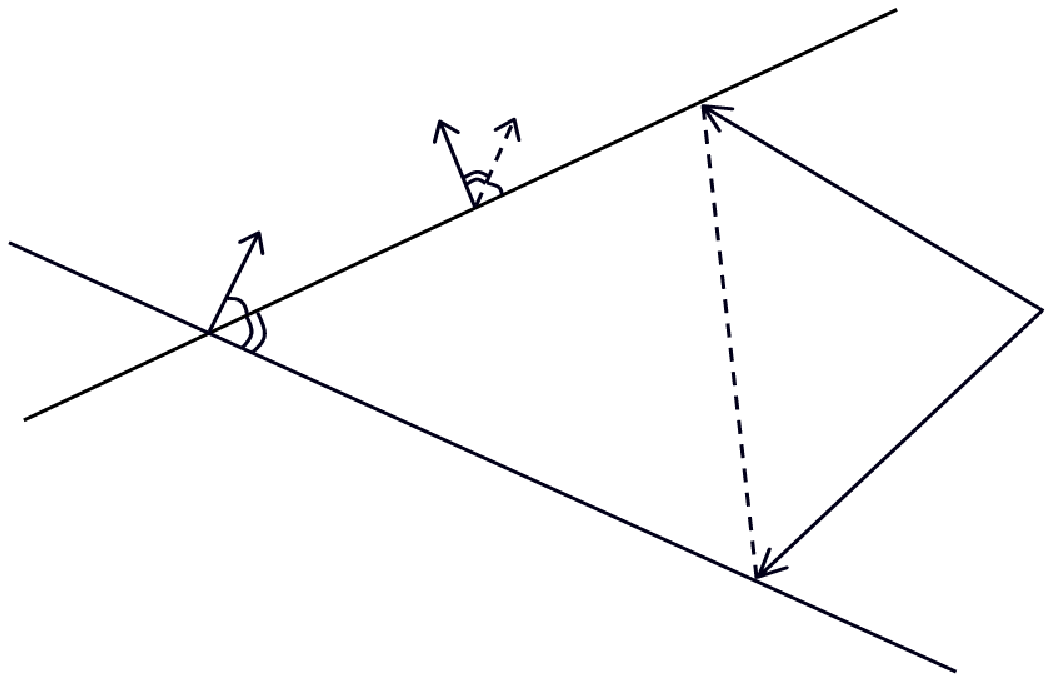}
\vskip -2.5cm \hskip 1.3cm $P(\hat{x},\hat{y})$ 
\vskip -0.9cm  \hskip 2.3cm $\theta_2 -\theta_1$ 
\vskip -0.5cm \hskip 7.8cm  $O(0,0)$
\vskip -2.3cm \hskip 4.0cm $Q_2(\alpha_2,\beta_2)$ 
\vskip 3.5cm \hskip 4.0cm $Q_1(\alpha_1,\beta_1)$
\vskip -5.0cm \hskip 7.0cm $\Gamma_{t_2}$
\vskip 5.cm \hskip 7.0cm $\Gamma_{t_1}$
\vskip -0.2cm
\centerline{Figure 1}
\bigskip

Observe that the first term of the right hand side in (2.14) is

\begin{align} 
-\frac{PQ_{1}sin{\psi}PQ_{2}sin{\psi}}{sin{\psi}}&=-(PQ_{1})(PQ_{2})sin{\psi}&\nonumber\\
&=-2Area({\Delta}PQ_{1}Q_{2}).&\nonumber
\end{align}

\noindent It is easy to see that the second term of the right hand side in (2.14) equals twice the negative area of the triangle $Q_{1}OQ_{2}$. Note that if $PQ_{1}OQ_{2}P$ is in the order of the clockwise direction, then the first and second terms of the right hand side in (2.14) are twice the positive area of ${\Delta}PQ_{1}Q_{2}$ and ${\Delta}Q_{1}OQ_{2}$, respectively. Without loss of generality, we may take $({\alpha},{\beta})$ to be proportional to $(cos{\theta},sin{\theta})$, the normal to $\Gamma_t$, as expressed in (2.7) with ${\xi}=0$. Then we can have a simpler expression for $\hat{z}_{2}-\hat{z}_{1}$: (denote ${\delta}(t_{i})$ by $\delta_{i}$)

\begin{align} 
&\hat{z}_{2}-\hat{z}_{1}=\frac{2{\delta_1}{\delta_2}-({\delta_1}^{2}+{\delta_2}^{2})cos(\theta_{2}-\theta_{1})}{sin(\theta_{2}-\theta_{1})}+{\gamma}_{2}-{\gamma}_{1}.&\tag{2.15}
\end{align}

\noindent Here we have combined the first and second terms of the right hand side in (2.14) using the equality $-(\alpha_{2}\beta_{1}-{\alpha}_{1}\beta_{2}) = {\delta_1}{\delta_2}sin(\theta_{2}-\theta_{1})$. The first term of the right hand side in (2.15) being the negative or positive twice area of the quadrilateral $PQ_{1}OQ_{2}$ (see Figure 2 below) can be deduced from elementary plane geometry as follows. 

\bigskip

\includegraphics[width=8.0cm]{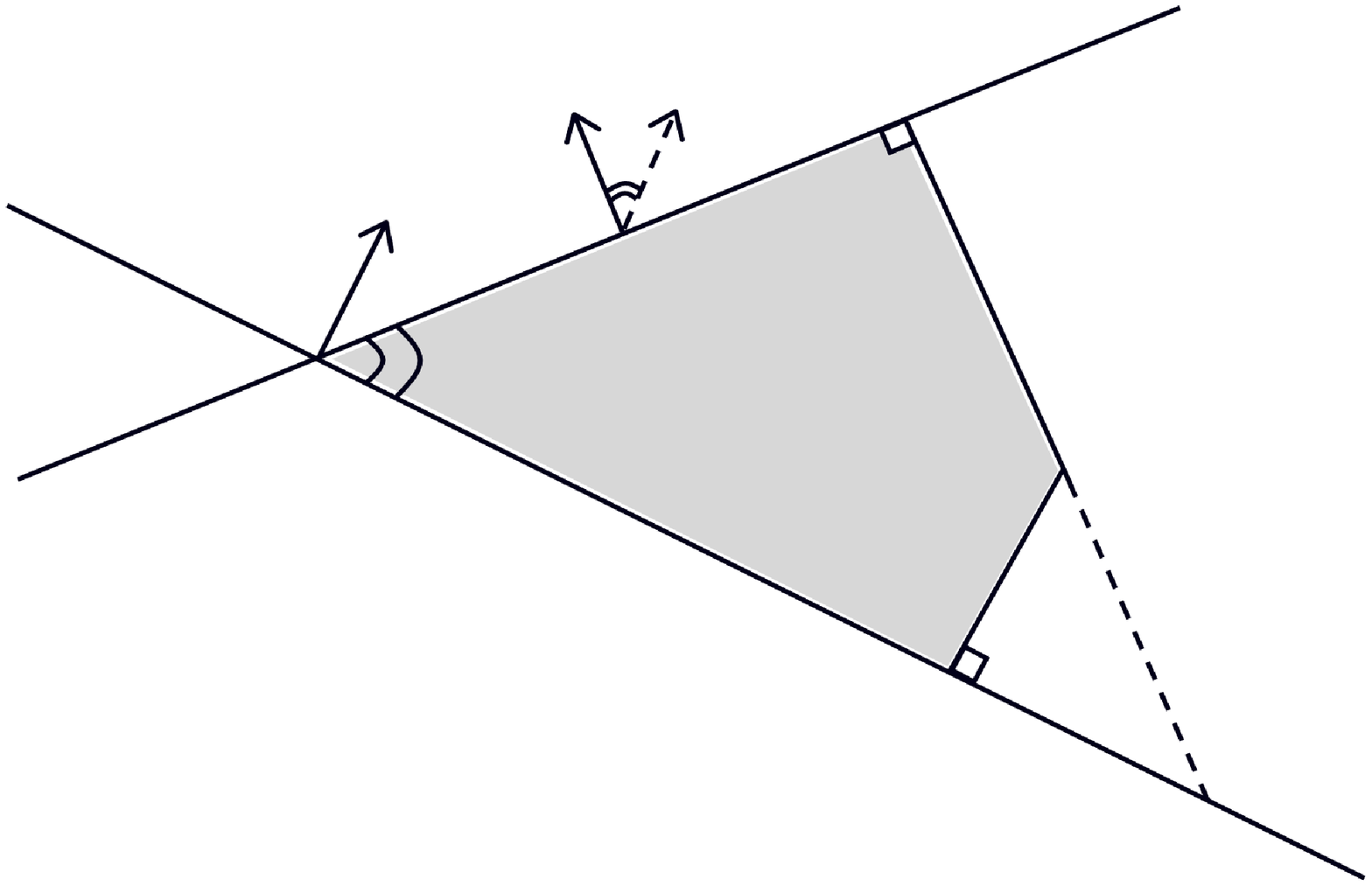}
\vskip -3.1cm \hskip 1.8cm $P$ 
\vskip -.8cm  \hskip 2.6cm $\phi$
\vskip -1.9cm \hskip 4.9cm $Q_2$
\vskip 3.1cm \hskip 5.3cm $Q_1$
\vskip -3.2cm \hskip 5.9cm $\delta_2>0$
\vskip 1.0cm \hskip 6.4cm $O(0,0)$
\vskip 0.5cm \hskip 5.9cm $\delta_1<0$
\vskip 0.9cm \hskip 6.9cm $R$
\vskip -5.4cm \hskip 7.0cm $\Gamma_{t_2}$
\vskip 4.9cm \hskip 7.4cm $\Gamma_{t_1}$
\vskip -0.2cm
\centerline{Figure 2}

\bigskip

Since $\Delta{PRQ_{2}}$ is similar to $\Delta{ORQ_{1}}$, we have $\frac{PQ_{2}}{|\delta_{1}|} = \frac{|\delta_{2}|+|\delta_{1}|(cos{\phi})^{-1}}{|\delta_{1}|tan{\phi}}$ ($\phi$ denotes the degree of the angle $\angle{Q_{2}PQ_{1}}$). Note that $\delta_{1},\delta_{2}$ may be negative. It follows that $PQ_{2} = (|\delta_{2}|+|\delta_{1}|(cos{\phi})^{-1})cot{\phi}$. Now we compute 

\begin{align}
2Area(PQ_{1}OQ_{2})&=2(Area(\Delta{PRQ_{2}})-Area(\Delta{ORQ_{1}}))&\tag{2.16}\\
&=(|\delta_{2}|+|\delta_{1}|(cos{\phi})^{-1})^{2}cot{\phi}-{\delta_{1}}^{2}tan{\phi}&\nonumber\\
&=\frac{2|\delta_{1}\delta_{2}|+({\delta_{1}}^{2}+{\delta_{2}}^{2})cos{\phi}}{sin{\phi}}.&\nonumber
\end{align}

\noindent For the situation shown in Figure 2, $PQ_{1}OQ_{2}P$ is in the order of the counterclockwise direction while $\delta_{1}$ is negative, $\delta_{2}$ is positive, and ${\phi}=\theta_{2}-\theta_{1}$. So $|\delta_{1}\delta_{2}|$ is in fact -$\delta_{1}\delta_{2}$.
 
\bigskip

{\bf Proof of part (c) of Theorem A:}

The map $X$ (given by (2.1), (2.2), and (2.3)) with the above choice of $\alpha$, $\beta$ is injective if and only if $\hat{z}_{2}-\hat{z}_{1}$ given by the formula (2.15) is not 0 in the case that $\Gamma_{t_1}\cap\Gamma_{t_2}$ consists of exactly one point $(\hat{x},\hat{y})$ for ${t_1}\neq{t_2}$. Suppose $\Gamma_{t_1}$ coincides with $\Gamma_{t_2}$ (in this case $sin(\theta_{2}-\theta_{1})=0$). Then $(\alpha_{1},\beta_{1}) = (\alpha_{2},\beta_{2})$. Therefore in this situation a necessary and sufficient condition for the injectivity of $X$ is $\gamma_{1}\neq\gamma_{2}$ by (2.11) (note that $\tilde{\Gamma}_{t_{1}}$ is parallel to, but not identified with, $\tilde{\Gamma}_{t_{2}}$ if $\gamma_{1}\neq\gamma_{2}$). We have proved part (c) of Theorem A.

\begin{flushright}
Q.E.D.
\end{flushright}

\medskip 

{\bf Example 2.1.} We take ${\alpha}=0$, ${\beta}(t)=t$, ${\gamma}(t)=-t$, and ${\theta}(t) = tan^{-1}t+\frac{\pi}{2}$ for $-\infty < t < +\infty$ in (2.1), (2.2), and (2.3). We claim that the associated $X=X(s,t)$ defined over $R{\times}R$ is a proper embedding. First observe that $cos{\theta}(t)=-\frac{t}{\sqrt{1+t^{2}}}$ and $sin{\theta}(t)=\frac{1}{\sqrt{1+t^{2}}}$. It is then easy to eliminate the parameters $s,t$ to get the equation in Cartesian coordinates as follows:

\begin{align}
&z(x+1)=y(x-1).&\tag{2.17}           
\end{align}

\noindent Now to prove that $X$ is an immersion, we argue that (2.8a) and (2.8b) can not hold simultaneously for some $(s,t)$ (note that $\delta = -cos\theta (t)=\frac{t}{\sqrt{1+t^{2}}}$ and $\xi = \frac{t^{2}}{\sqrt{1+t^{2}}}$). Since ${\gamma}(t)=-t$ is monotonically decreasing and the first two terms of the right hand side in (2.14) are both nonpositive (-2$Area({PQ_{1}OQ_{2}})$), we conclude that $\hat{z}_{2}-\hat{z}_{1}$ is always negative for $t_{2}>t_{1}$. By part (c) of Theorem A (another situation never occurs), we know that $X$ is injective. Observe that $|X(s,t)| \longrightarrow \infty$ as $|(s,t)| \longrightarrow \infty$. We can then easily show that $X$ is proper and a homeomorphism between $R{\times}R$ and the image of $X$ with the topology induced from $R^3$. Thus $X$ is a proper embedding. 

Also this is an example of helicoid type. Namely, the characteristic lines $\tilde{\Gamma}_{t}$ are lying in a family of parallel planes as explained below. Observe that $(sin{\theta}(t),-cos{\theta}(t),$ ${\beta}(t)sin{\theta}(t)+{\alpha}(t)cos{\theta}(t))$ is tangent to $\tilde{\Gamma}_{t}$. The cross product of two such vectors must be parallel to a constant vector for all $t_{1}\neq{t_2}$ if $\tilde{\Gamma}_{t}$ are lying in a family of parallel planes. Computing the cross product of such vectors for $t_{1}=0$ and $t_{2}\neq{0}$, we obtain the vector $(0, -t_{2}sin{\theta}(t_{2}), -cos{\theta}(t_{2}))$. This is parallel to a constant vector $(0, -1, 1)$ for any $t_2$. In Example 3.1, we will express this helicoid type p-minimal surface in the form of (1.2). 

\medskip

We now give a direct application of (2.15) and (2.16). Denote $\Gamma_{t_i}$ by $\Gamma_{i}$ for $1{\leq}i{\leq}3$. Suppose $\Gamma_{1}{\cap}\Gamma_{2}=\{P\}$, $\Gamma_{2}{\cap}\Gamma_{3}=\{Q\}$, and $\Gamma_{3}{\cap}\Gamma_{1}=\{R\}$ (see Figure 3). 

\bigskip

\includegraphics[width=8.0cm]{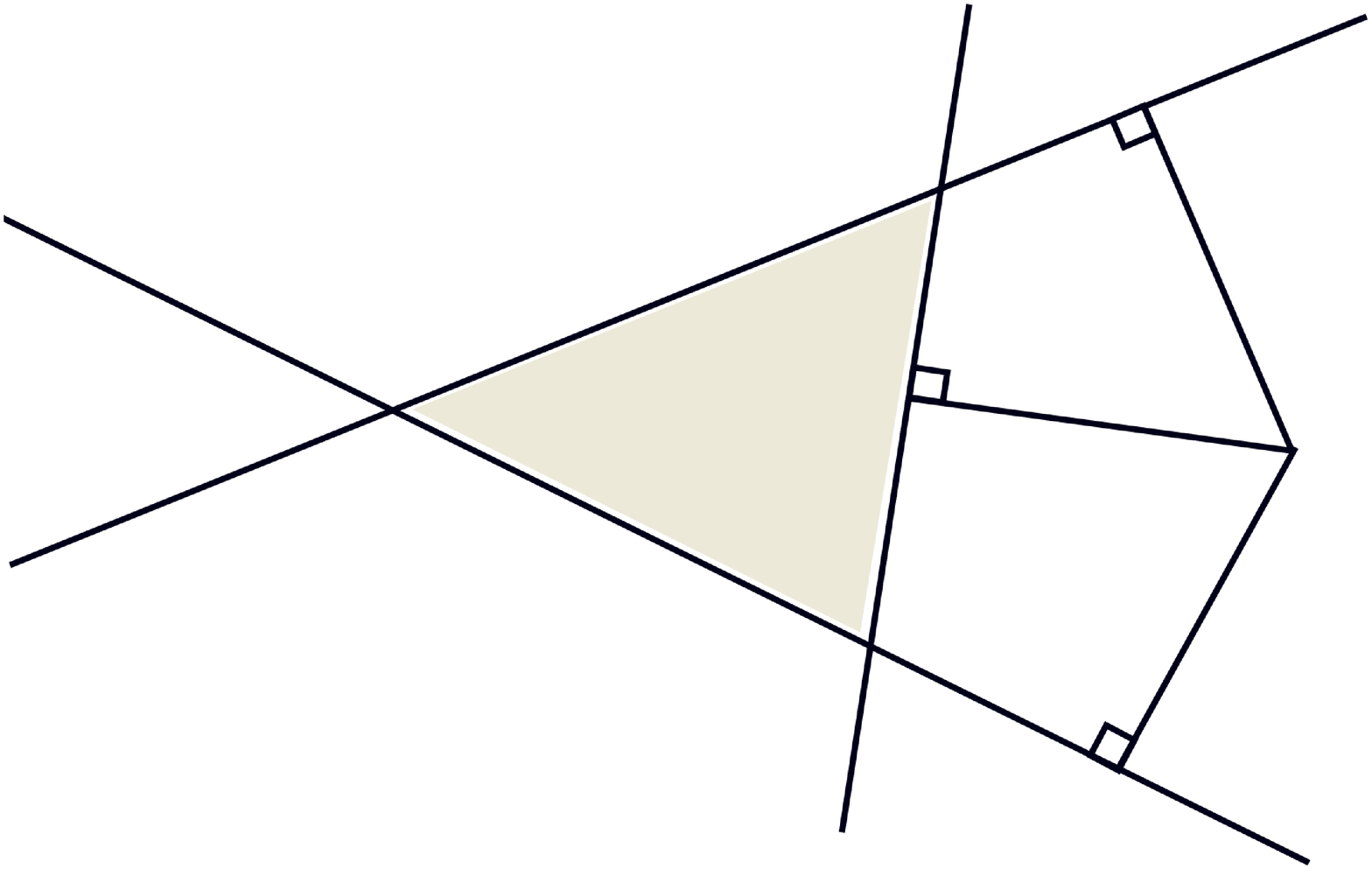}
\vskip -3.0cm \hskip 2.0cm $P$ 
\vskip -2.0cm \hskip 4.9cm $Q$
\vskip 1.1cm \hskip 5.4cm $U$
\vskip 1.1cm \hskip 5.0cm $R$
\vskip -4.7cm \hskip 5.3cm $\Gamma_3$
\vskip 0.1cm \hskip 6.5cm $S$
\vskip 1.7cm \hskip 7.5cm $O(0,0)$
\vskip 1.7cm \hskip 6.2cm $T$
\vskip -5.0cm \hskip 7.0cm $\Gamma_2$
\vskip 4.9cm \hskip 7.4cm $\Gamma_1$
\vskip -0.2cm
\centerline{Figure 3}

\bigskip

\noindent {\bf Proposition 2.2.} {\it Suppose we have the situation as described above. Then the sum of the height differences at the intersection points is twice the area of the triangle region surrounded by $\Gamma_1$, $\Gamma_2$, and $\Gamma_3$. That is to say,}

\begin{align}
&(z_{1}-z_{2})(P)+(z_{2}-z_{3})(Q)+(z_{3}-z_{1})(R)=2Area\Delta{PRQ}.&\tag{2.18}
\end{align}

\bigskip

{\bf Proof of Proposition 2.2:}

Without loss of generality, we consider the situation as shown in Figure 3. By (2.15) and (2.16) we obtain that $(z_{2}-z_{1})(P) = -2Area(PTOS) + \gamma_{2}-\gamma_{1}$ since $PTOS$ is in the order of the counterclockwise direction (otherwise we would have the positive sign in front of $2Area(PTOS)$). Similarly we have $(z_{3}-z_{2})(Q) = +2Area(QUOS) + \gamma_{3}-\gamma_{2}$ and $(z_{1}-z_{3})(R) = +2Area(RTOU) + \gamma_{1}-\gamma_{3}$. Summing up these three identities gives (2.18).

We can also prove (2.18) by integrating the contact form (restricted to a local graph) along the boundary of $\Delta{PRQ}$. Recall that the contact form ${\Theta} = dz+xdy-ydx$ ([CHMY]) vanishes along a characteristic curve. It follows that

\begin{align}
0&=z_{1}(R)-z_{1}(P)+z_{3}(Q)-z_{3}(R)+z_{2}(P)-z_{2}(Q)+\oint_
{{\partial}({\Delta}PRQ)}(xdy-ydx)&\nonumber\\ 
&=(z_{2}-z_{1})(P)+(z_{3}-z_{2})(Q)+(z_{1}-z_{3})(R)+2Area{\Delta}PRQ&\nonumber
\end{align}

\noindent where we have used Green's or Stokes' theorem. We are done.

\begin{flushright}
Q.E.D.
\end{flushright} 

\bigskip

\section{{\bf Helicoid type p-minimal surfaces: Proofs of Theorem B and Theorem C}}

We first prove Theorem B. There are two cases:

Case 1. $P_t$'s are not perpendicular to the $xy$-plane.

\noindent We observe that the Legendrian line ${\tilde{\Gamma}}_{t}$ always contains a point $p_{(t)}$ such that $P_t$ is the contact plane past $p_{(t)}$. So if we write $p_{(t)} = (x_{0}(t),y_{0}(t),z_{0}(t))$, the vector $(-y_{0}(t),x_{0}(t),1)$ is normal to $P_t$. Since $P_t$'s are parallel to each other, we conclude that $(x_{0}(t),y_{0}(t))$ is a constant vector $(x_{0},y_{0})$ for all $t$. Without loss of generality, we may assume that $z_{0}(t)=t$ and hence $P_t$ is defined by $z = y_{0}x-x_{0}y+t$. Therefore we can describe $\Sigma$ in parameters $s$ and $t$ as follows (see (2.1), (2.2), (2.3) in Section 2 or (4.9) in [CHMY]):

\begin{align}
&x=s(sin{\theta}(t))+x_{0}&\tag{3.1}\\
&y=-s(cos{\theta}(t))+y_{0}& \tag{3.2}\\
&z=s[y_{0}sin{\theta}(t)+x_{0}cos{\theta}(t)]+t& \tag{3.3}
\end{align}

\noindent where ${\theta}\in C^2(R)$. From (3.1) and (3.2) we obtain

\begin{align}
&(x-x_{0})cos{\theta}(t)+(y-y_{0})sin{\theta}(t)=0.&\tag{3.4}
\end{align}

\noindent Substituting $s(sin{\theta}(t))=x-x_{0}$ and $s(cos{\theta}(t))=-(y-y_{0})$ into (3.3) gives

\begin{align}
&t=z-y_{0}(x-x_{0})+x_{0}(y-y_{0})=z-y_{0}x+x_{0}y.&\tag{3.5}
\end{align}

\noindent We have shown that $\Sigma$ is of the form (1.2) by (3.4) and (3.5).

Case 2. $P_t$'s are perpendicular to the $xy$-plane.

We may assume $P_t$'s are defined by $-bx+ay = t$ with $a^{2}+b^{2} = 1$ without loss of generality. So we can describe $\Sigma$ in parameters $s$ and $t$ as follows:

\begin{align}
&x=as-bt&\tag{3.6}\\
&y=bs+at& \tag{3.7}\\
&z=f(t)s+g(t).& \tag{3.8}
\end{align}

\noindent ${\tilde{\Gamma}}_{t}$ being Legendrian (i.e., dz+xdy-ydx=0 along ${\tilde{\Gamma}}_{t}$) implies that $f(t)=t$. Expressing $s$, $t$ in $x$, $y$ by (3.6),(3.7) and then substituting the result into (3.8), we obtain (1.1). We have finished the proof of Theorem B.

Next we will give a proof of Theorem C. Let
$F(x,y,z)=(x-x_{0})cos{\theta}(t)+(y-y_{0})sin{\theta}(t)$ be a defining
function with $t$ given by (3.5) for $\Sigma$. Compute the gradient $\nabla F$
of $F$. With the substitution of $x-x_{0}=s(sin{\theta}(t))$ and
$y-y_{0}=-s(cos{\theta}(t))$, we obtain

\begin{align}
&{\nabla F}=(sy_{0}{\theta}^{\prime}(t)+cos{\theta}(t),
-sx_{0}{\theta}^{\prime}(t)+sin{\theta}(t), -s{\theta}^{\prime}(t)).&\tag{3.9}
\end{align}

\noindent It is easy to see that ${\nabla F}$ never vanishes. Moreover, a point
$(x,y,z)\in \Sigma$ is a singular point if and only if ${\nabla}_{b}F{\equiv}(\hat{e}_{1}F)\hat{e}_{1}+(\hat{e}_{2}F)\hat{e}_{2}=0$ where 
$\hat{e}_{1}=\frac{\partial}{\partial{x}}+y\frac{\partial}{\partial{z}},\hat{e}_{2}
=\frac{\partial}{\partial{y}}-x\frac{\partial}{\partial{z}}$. So the condition for 
$(x,y,z)$ to be a singular point becomes the following one by
(3.9) through a straightforward computation:

\begin{align}
&cos{\theta}(t)(1+s^{2}{\theta}^{\prime}(t))=0,
sin{\theta}(t)(1+s^{2}{\theta}^{\prime}(t))=0.&\tag{3.10}
\end{align}

Therefore if ${\theta}^{\prime}(t)\geq 0$, then
$1+s^{2}{\theta}^{\prime}(t)\geq 1$. It follows from (3.10) that
$cos{\theta}(t)=sin{\theta}(t)=0$, a contradiction. Hence $\Sigma$ has no
singular point. Conversely, if ${\theta}^{\prime}(t)<0$ for some $t$, we can
take $s^2=-\frac{1}{{\theta}^{\prime}(t)}$ such that
$1+s^{2}{\theta}^{\prime}(t)=0$. So (3.10) holds, and hence $\Sigma$ has a
singular point. We have proved Theorem C.

\bigskip

{\bf Example 3.1.} In Example 2.1, we learned that the surface $z(x+1)=y(x-1)$ (see (2.17)) is a $C^2$ smooth, helicoid type p-minimal surface. So by Theorem B, we should be able to express it in the form of (1.2). Since $(0,-1,1)$ is normal to the associated parallel planes (see Example 2.1), we know that $(x_0,y_0)=(-1,0)$ from the proof of Theorem B case 1. Now from (1.2) we find that the associated $\theta (t)$ must satisfy $tan\theta (t)=-\frac{x+1}{y}$. On the other hand, $t=z-y=-\frac{2y}{x+1}$ by (2.17). It follows that $tan\theta (t)=\frac{2}{t}$, hence $\theta (t)=cot^{-1}(\frac{t}{2})$.   

\bigskip

\section{{\bf Properness and non-helicoid type p-minimal surfaces}}

Let $X:\Omega\rightarrow H_{1}$ be a $C^2$ smooth immersion of the p-minimal surface $\Sigma =X(\Omega )$ in $H_1$. Here the domain $\Omega$ denotes either $R^{2}=R^{1}\times R^{1}$ (band type) or $R^{1}\times S^{1}$ (annulus type) and the first parameter describes the rulings of $\Sigma$. We say that $\Sigma$ is proper if $X$ is proper, i.e., $X^{-1}(K\cap\Sigma)$ is compact for any compact subset $K$ of $H_1$. We say that $\Sigma$ is complete if it is complete with respect to the metric induced from the standard Euclidean metric of $R^3$ identified with $H_1$. Let $\tilde{\Gamma}_{t}'s$ denote the rulings of $\Sigma$ ($t$ is the second parameter). Denote the distance between $\tilde{\Gamma}_{t}$ and the origin $O$ by $r_t$. We have the following result.

\bigskip

\noindent {\bf Proposition 4.1.} {\it (1) Suppose that an annulus type (i.e., $\Omega =R^1\times S^1$) immersed p-minimal surface $\Sigma$ is complete. Then $\Sigma$ is proper.}

{\it (2) Suppose that a band type (i.e., $\Omega =R^1\times R^1$) immersed p-minimal surface $\Sigma$ is complete. Then $\Sigma$ is proper if and only if  $lim_{t\rightarrow\pm\infty}r_{t}=+\infty$.}

\bigskip

{\bf Proof: } First observe that $\Sigma$ is proper if and only if $X^{-1}(\bar{B}(p,\rho )\cap\Sigma)$ is compact (hence $\bar{B}(p,\rho )\cap\Sigma$ is also compact) for any closed ball $\bar{B}(p,\rho )$ of center $p$ and radius $\rho$ in $R^3$. The completeness of $\Sigma$ implies that $\tilde{\Gamma}_{t}\subset\Sigma$ must be the whole straight line.

For the proof of (2), if $\Sigma$ is proper, then the compactness of $X^{-1}(\bar{B}(O,\rho )\cap\Sigma)$ in $R^2$ implies its boundedness for any large $\rho$. It follows that $\tilde{\Gamma}_{t}$ lies outside the closed ball $\bar{B}(O,\rho )$ for $|t|$ large enough. This is just what $lim_{t\rightarrow\pm\infty}r_{t}=+\infty$ means. Conversely, the condition $lim_{t\rightarrow\pm\infty}r_{t}=+\infty$ implies that for any $p$, $\rho >0$, there exist $t_{1},t_{2}\in R$, $t_{1}<t_{2}$, such that

\begin{align}
&\bar{B}(p,\rho )\cap\Sigma\subset\bigcup_{t\in [t_{1},t_{2}]}\tilde{\Gamma}_{t}.&\tag{4.1}
\end{align}

\noindent Since $\bar{B}(p,\rho )\cap\tilde{\Gamma}_{t}$ is a compact subset of $\tilde{\Gamma}_{t}$, it is not hard to show that $X^{-1}(\bar{B}(p,\rho )\cap\Sigma)\subset [s_{1},s_{2}]\times [t_{1},t_{2}]$. On the other hand, we can easily see that $\bar{B}(p,\rho )\cap\Sigma$ is a closed subset of $H_1$ due to the completeness of $\Sigma$. So $X^{-1}(\bar{B}(p,\rho )\cap\Sigma)$ is a closed subset of $R^2$. We have shown that $X^{-1}(\bar{B}(p,\rho )\cap\Sigma)$ is bounded and closed in $R^2$, hence compact. 

For the proof of (1), we note that (4.1) holds true with $[t_{1},t_{2}]$ replaced by $S^1$ without any condition (the right hand side is the whole surface $\Sigma$). A similar argument as in the above proof of (2) shows that $X^{-1}(\bar{B}(p,\rho )\cap\Sigma)$ is bounded in $R^{2}\supset R^1\times S^1$ and closed in $R^1\times S^1$, hence closed in $R^{2}$.

\begin{flushright}
Q.E.D.
\end{flushright}      

\bigskip

In the remaining section, we will discuss the situation of a p-minimal surface being not of helicoid type. The first possibility is that the union of the parallel planes containing characteristic lines is no longer the whole space. This can occur as shown in the following example. 

\medskip

{\bf Example 4.1.} We consider the surface $\Sigma$ defined by $(z-xy)^{2}=y$ in $R^3$. This surface can be viewed as an example of $(1.1)$ with $a=1$, $b=0$, and a two-valued function $g=\pm\sqrt{y}$. We can easily verify that $\Sigma$ is a $C^2$ smooth, properly embedded p-minimal surface with the characteristic lines $\{ y=t^{2},z=t^{2}x+t,-\infty <t<+\infty\}$ lying in the parallel planes $P_{t}=\{y=t^{2}\}$. Note that $P_{t}=P_{-t}$ and $\bigcup_{t\in R}P_{t}=\{ y\geq 0\}\neq R^{3}$ or $H_1$. Since the characteristic lines are still lying in the parallel planes, such a surface is called of helicoid type in the weak sense. 

\medskip

{\bf Example 4.2.} Take $(\alpha (t),\beta (t))=(cos\theta (t),sin\theta (t))$ where $max\theta (t) - min\theta (t) < \pi$ in (2.7) (i.e. $\delta =1,\xi =0$) and $\gamma (t)$ such that $lim_{t\rightarrow\pm\infty}|\gamma (t)|=+\infty$, ${\gamma}^{\prime}(t)$ has the same sign as ${\theta}^{\prime}(t)$, and ${\theta}^{\prime}(t)$ never vanishes (for instance, $\theta (t)=tan^{-1}t$ and $\gamma (t)=t$). Namely, we consider the $C^2$ smooth map $X:R^{1}\times R^{1}\rightarrow H_{1}$ defined by $X(s,t) = (x(s,t),y(s,t),z(s,t))$ where $x,y,z$ are given in (2.1), (2.2), (2.3) with the above data (assuming $\theta ,\gamma \in C^{2}$). Observe that $z=s+\gamma (t)$. It follows that 

$$x^{2}+y^{2}+z^{2}=1+s^{2}+(s+\gamma (t))^{2}$$.

\noindent Therefore $r_{t}=1+\frac{1}{2}{\gamma}^{2}(t)$ by elementary calculus, and hence $lim_{t\rightarrow\pm\infty}r_{t}=+\infty$ by our assumption on $\gamma$. Now we can show that $X$ is proper by a similar argument as in the proof of Proposition 4.1 (noting that we only need 
that $\tilde{\Gamma}_{t}\subset\Sigma =X(R^{1}\times R^{1})$ is a whole straight line instead of $\Sigma$ being complete). We can examine the conditions in Theorem A (a) and (c) respectively to conclude that $X$ is immersed and injective. A proper injective immersion must be an embedding. 

Next we claim that the characteristic lines $\tilde{\Gamma}_{t}'s$ are not lying in the parallel planes $P_{t}'s$. If they are, we follow the argument as in the proof of Theorem B to reach a contradiction. First Case 2 is ruled out apparently. In Case 1, $\tilde{\Gamma}_{t}$ must contain a point having constant $x,y$ components for any $t$. This is clearly impossible since $\Gamma_t$'s are not having a common intersection point. We have constructed a family of $C^2$ smooth, band type, properly embedded p-minimal surfaces in $H_1$, which are not of helicoid type even in the weak sense (i.e., the characteristic lines are not lying in the parallel planes).

The condition (2.9) for a singular point reduces to $(s^{2}+1){\theta}^{\prime}+{\gamma}^{\prime}=0$. Since ${\theta}^{\prime}$ and ${\gamma}^{\prime}$ have the same sign by assumption, the surfaces under consideration are all having no singular points.

\medskip

Example 4.1 shows that a properly embedded p-minimal surface can exist in the "right" halfspace. However, this can not occur for the upper halfspace (see Theorem D).
      
\bigskip

{\bf Proof of Theorem D:}

The completeness implies that $\Sigma =\bigcup \tilde{\Gamma}_t$ where each $\tilde{\Gamma}_t$ is a whole characteristic line sitting in $P_{+}$ (or $P_{-}$), hence parallel to the contact plane $P$. Let $P_{t}$ denote the contact plane parallel to $P$, containing $\tilde{\Gamma}_{t}$. Note that $P_{t}$'s are not perpendicular to the $xy$-plane since $\partial /\partial z$ is not annihilated by the contact form $dz+xdy-ydx$. So a similar argument in the proof of Theorem B, Case 1 shows that $\tilde{\Gamma}_{t}$ must pass through $(x_{0},y_{0},\gamma (t))$. Let $z_{0}$ denote the minimum (or the maximum if in $P_{-}$) of $\gamma (t)$ over $t$. Our assumptions on $\Sigma$ force that $\gamma (t)=z_{0}$ for all $t$, and $\Sigma$ is the contact plane past $(x_{0},y_{0},z_{0})$, parallel to $P$.   

\begin{flushright}
Q.E.D.
\end{flushright}

\bigskip







\begin{thebibliography}{W}

\bibitem[CHMY]{chmy}
Cheng, J.-H., Hwang, J.-F., Malchiodi, A., and Yang, P., {\em Minimal surfaces
in pseudohermitian geometry}, arXiv: math.DG/0401136.

\bibitem[Col]{col}
Collin, P., {\em Topologie et courbure des surfaces minimales de $R^3$}, Annals
of Math., 2nd Series, 145 (1997) 1-31.

\bibitem[Cos]{cos}
Costa, C., {\em Example of a complete minimal immersion in $R^3$ of genus one and
three embedded ends}, Bull. Soc. Bras. Mat., 15 (1984) 47-54.

\bibitem[DHKW]{dhkw}
Dierkes, U., Hildebrandt, S., K$\ddot{u}$ster, A., and Wohlrab, O., {\em Minimal 
surfaces I, II}, G.M.W., 295, 296, Springer-Verlag, 1992.

\bibitem[Fan]{fan}
Fang, Y., {\em Lectures on minimal surfaces in $R^3$}, Proceedings of CMA, the 
Australian National University, Vol. 35, 1996.

\bibitem[GP1]{gp1}
Garofalo, N. and Pauls, S., {\em The Bernstein problem in the Heisenberg group}, 
arXiv:math.DG/0209065.

\bibitem[GP2]{gp2}
Garofalo, N. and Pauls, S., {\em The Bernstein problem in the Heisenberg group}, 
preprint, 2004.


\bibitem[HKW]{hkw}
Hoffman, D., Karcher, H., and Wei, F., {\em The genus one helicoid and the minimal
surfaces that led to its discovery}, In {\it Global Analysis and Modern Mathematics},
Publish or Perish Press, ed. K. Uhlenbeck (1993) 119-170.

\bibitem[HM1]{hm1}
Hoffman, D. and Meeks III, W. H., {\em Complete embedded minimal surfaces of finite
total curvature}, Bull. A.M.S., 12 (1985) 134-136.

\bibitem[HM2]{hm2}
Hoffman, D. and Meeks III, W. H., {\em The strong halfspace theorem for minimal surfaces}, Inventiones Math.,101 (1990) 373-377.

\bibitem[Mon]{mon}
Monge, G., {\em Application de l'analyse \`{a} la g\'{e}om\'{e}trie}, Paris,
Bachelier, 1850.

\bibitem[MR]{mr}
Meeks III, W. H. and Rosenberg, H., {\em The uniqueness of the helicoid and the 
asymptotic geometry of properly embedded minimal surfaces with finite topology}, preprint.

\bibitem[Nit]{nit}
Nitsche, J. C. C.,
{\it Vorlesungen $\ddot{u}$ber Minimalfl$\ddot{a}$chen}, Grundelehren der Math., Bd. 199, Springer- Verlag, Berlin- New York, 1975. 

\bibitem[Oss]{oss}
Osserman, R.,
{\it A survey of minimal surfaces}, Dover Publications, Inc., New York, 1986.

\bibitem[Pau]{pau}
Pauls, S. D., {\em Minimal surfaces in the Heisenberg group}, Geometric Dedicata, 104
(2004) 201-231. 

\end{thebibliography}
\end{document}